\documentclass[12pt]{amsart}
\usepackage{epsfig,color}

\makeatother

\theoremstyle{definition}

\def\fnum{equation} 
\newtheorem{Thm}[\fnum]{Theorem}

\numberwithin{equation}{section}

\newcommand{\R}{{\text{R}}}

\newcommand{\Ric}{{\text{Ric}}}

\newcommand{\Tr}{{\text{Tr}}}

\newcommand{\Hess}{{\text {Hess}}}

\def\RR{{\bold R}}
\def\RP{{\bold{RP}}}
\def\SS{{\bold S}}

\newcommand{\dv}{{\text {div}}}
\newcommand{\e}{{\text {e}}}

\newcommand{\cL}{{\mathcal{L}}}

\newcommand{\cP}{{\mathcal{P}}}

\newcommand{\eqr}[1]{(\ref{#1})}

\title{Singularities   and diffeomorphisms}

\author{Tobias Holck Colding}%
\address{MIT, Dept. of Math.\\
77 Massachusetts Avenue, Cambridge, MA 02139-4307.}
\author{William P. Minicozzi II}%

\thanks{The  authors
were partially supported by NSF  DMS Grants 2104349 and 2005345.}

\email{colding@math.mit.edu and minicozz@math.mit.edu}

\begin{document}

\maketitle

{\centering\footnotesize Dedicated to Blaine Lawson with admiration.\par}

 \begin{abstract}
   Comparing and recognizing metrics can be extraordinarily difficult
because of the group of diffeomorphisms.  Two metrics, that could even be the same, could look completely different in different coordinates.  This is the gauge problem.  The general gauge problem is extremely subtle
for non-compact spaces. Often it can be avoided if one uses some additional structure of the particular situation.  However, in many  problems there is no additional structure.  Instead we solve the gauge problem directly in great generality.

   The techniques and ideas apply to many problems.   We use them to solve a well-known open problem in Ricci flow.   
 
We solve the gauge problem by solving a nonlinear system of PDEs.  The PDE produces a diffeomorphism that fixes an appropriate gauge in the spirit of the slice theorem for group actions.  We then show optimal bounds for the displacement function of the diffeomorphism.  
  \end{abstract}

\section{Introduction}

 Suppose we have two weighted manifolds $(M_i,g_i,f_i)$ for $i=1$, $2$ satisfying some PDE.
 Assume that on a large, but compact set, the manifolds $M_i$, metrics $g_i$ and weights $\e^{-f_i}$ almost agree after identification by a diffeomorphism.  
\begin{itemize}
\item Are the manifolds, metrics and weights the same everywhere after some identification?
\end{itemize}
This is a common problem in many questions.
  The major obstacle for understanding this in general 
  is the infinite dimensional gauge group of diffeomorphisms:
\begin{itemize}
\item Two metrics, that could even be the same, could look very different in different coordinates.
\end{itemize}

In some situations the gauge problem can be avoided if there is some additional structure.  A classical example 
is the  Killing-Hopf theorem that classifies constant curvature metrics.   This classification uses that the curvature tensor is constant to construct a ``canonical'' isometry between the two spaces.
In general, the gauge problem can be solved when there is strong asymptotic decay and circumvented when the space is characterized in a coordinate-free way, such as a large symmetry group,  the vanishing of a special tensor, or  a strong curvature condition.

In the problems we will be interested in,  the manifold will be non-compact and we will not have any special structure.
Thus, we will be forced to deal with the gauge problem head on.  
The flip side of this is that once we do that it gives new tools with broad applications.

\subsection{Where do  questions like these arise?}  

Problems about identifying spaces occur in many different situations.  The one we are interested in 
here comes from Ricci flow.  
A one parameter family $(M,g(t))$ of manifolds flows by the Ricci flow, \cite{H}, if 
$$g_t=-2\,\Ric_{g(t)}\,  ,$$
where $\Ric_{g(t)}$ is the Ricci curvature of the evolving metric $g(t)$ and $g_t$ is the time derivative of the metric.

The key to understand Ricci flow is the singularities that form. The simplest singularity is a homothetically shrinking sphere that disappears (beomes extinct) at a point.  The product of a sphere with $\RR$ gives a shrinking cylinder; this singularity, called a neck pinch, is much more complicated than the spherical extinction.  In dimension three, spherical extinctions and neck pinches are essentially the only singularities.  
Adding another $\RR$ factor, gives the so-called bubble sheet singularity that is only recently partially understood.  With each additional $\RR$ factor, the singularities become more complicated and the sets where they occur are larger.

A triple $(M,g,f)$ of a manifold $M$, metric $g$ and function $f$ on $M$ 
is a gradient shrinking Ricci soliton (or {\it shrinker}) if it satisfies
$$
\Ric + \Hess_f = \frac{1}{2} \, g \, .
$$
  Shrinkers give special solutions of the Ricci flow that evolve by rescaling up to diffeomorphism and are  singularity models.  They arise as time-slices of limits of rescalings (magnifications) of the flow around a fixed future singular point in space-time.   Such limits are said to be  tangent flows at the singularity.  Even when $M$ is compact, the shrinker is typically non-compact and the convergence is on compact subsets.
  Shrinkers also arise in other important ways, such as  blowdowns from $-\infty$ for ancient flows with bounded entropy.   Ancient flows are flows that have existed for all prior times; every tangent flow is ancient.  Shrinkers are the key singularities in Ricci flow and will be our focus here.

  Among shrinkers,  cylinders are particularly important; they are the most prevalent.    This is because the
   Almgren-Federer-White dimension reduction divides the singular set into strata whose dimension is the dimension of translation-invariance of the blowup.
   For Ricci flow, this
   suggests:
\begin{itemize}
\item Top strata of the singular set corresponds to points where the blowup is $\RR^{n-2} \times N^2$.
\item The next strata consists of points where the blowup is $\RR^{n-3} \times N^3$.
\end{itemize}
The $N$'s  are themselves shrinkers and have been classified in low dimensions by Cao-Chen-Zhu, Hamilton, Ivey, 
Naber, Ni-Wallach, Perelman.    
In dimensions two and three, they are $N^2 = \SS^2$ or $\RP^2$ and $N^3 = \SS^3$ or $\SS^2 \times \RR$ (plus quotients).
 The classification in dimension three
relies on an equation for the $2$-tensor $\frac{\Ric}{S}$ that fails in higher dimensions where there is no similar classification.  In fact, there are
huge families of shrinkers in higher dimensions.
Combining  dimension reduction with the classification in low dimensions, we see that
the most prevalent singularities are:\\
\centerline{$\SS^2\times \RR^{n-2}$ followed by $\SS^3\times \RR^{n-3}$ (and quotients).}

As one approaches a singularity in the flow and magnifies, one would like to know which singularity it is.  Since most singularities are non-compact yet the evolving manifolds are closed, one only sees a compact piece of the singularity at each time as one approaches it.  The next theorem recognizes the most prevalent singularities from just a compact piece.

\begin{Thm}	\label{t:rigid}
\cite{CM1} 
Cylindrical shrinkers  $\SS^{\ell} \times \RR^{n-\ell}$ are strongly rigid for any $\ell$.
\end{Thm}

The theorem holds for products of $\RR^{n-\ell}$ with quotients of $\SS^{\ell}$ and a large class of other 
positive Einstein manifolds; see \cite{CM1} for details.
{\it Strong rigidity means} that if another shrinker is close enough on a large compact set, then it must agree.  

In most problems in geometric PDEs, it would be unthinkable   to control an entire solution from just knowing roughly how it looks on a compact set.  This is exactly what we do here.  If one knew exactly how it looked like on a compact set, it would be much less surprising and essentially follow from unique continuation.  The surprising thing here is that we only assume closeness and only on a compact set and this enough to characterize the shrinker.  This is an illustration of a 
{\it{shrinker principle}}  which roughly says that  ``uniqueness radiates outwards''.  Nothing like this is true for   Einstein manifolds (or steady solitons), where gravitational instantons contain arbitrarily large arbitrarily Euclidean regions.
 The  shrinker principle was originally discovered in mean curvature flow 
  \cite{CIM}, \cite{CM2}. It has been conjectured since that something similar  holds for Ricci flow, but the gauge group has been one of the major obstacles.  Using extrinsic coordinates, the gauge is circumvented in mean curvature flow.

Tangent flows are  limits of a subsequence of rescalings at the singularity.  A priori different subsequences 
might give different limits.
Using Theorem \ref{t:rigid}, we get the following uniqueness:

\begin{Thm}	\label{t:uni}
\cite{CM1} 
For a Ricci flow, if  one tangent flow at a point in space-time is a   cylinder, then all other tangent flows at that point are also cylinders.
\end{Thm}

Unlike most results in Ricci flow, these results hold  for every $n$ and $\ell$. Increasing the dimension of the Euclidean factor is a subtle problem (e.g.~surgery, cylindrical estimates, and  $k$-convexity estimates only allow small Euclidean factors).
For general $n$ and $\ell$, cylinders do not have a coordinate-free characterization.
  This is a major part of the difficulty, forcing us to address the gauge problem head on.

At singularities where the tangent flows are compact shrinkers the singularities are isolated in space-time.  For compact shrinkers, rigidity was proven in dimension three by Hamilton in 1982 and by Huisken in 1985 for higher dimensional spheres.

Rigidity for necks  $\SS^{n-1} \times \RR$ was proven independently by Li-Wang \cite{LW}.  They
are able to circumvent the gauge problem by using
 that  their Euclidean factor is a line.  They do that, in part, by using 
 tensors with  special properties on the product of a sphere with a  line to prove asymptotic structure and approximate symmetry.  Once they have this, 
 they are able to use again that their Euclidean factor is a line to 
 apply   Brendle's symmetry improvement to get $O(n)$ symmetry and, finally,
  Kotschwar's classification of rotationally symmetric shrinkers. 
  
  \subsection{Further applications}  Rigidity and uniqueness of blowups are fundamental questions in regularity theory that  
 have many  applications.    For instance in mean curvature flow, they play a major role in understanding the singular set, proving optimal regularity, 
 understanding solitons, 
 classifying ancient solutions, and understanding low entropy flows.

\section{What is needed for rigidity?}
We need to show that if two shrinkers are close on a large but compact set, then they agree identically everywhere.  This will be done by iterating two estimates: extension and improvement.  Extension  shows that the shrinkers remain close even on a fixed larger scale, but with a loss in the estimates.  Improvement recovers this loss and shows that they are even closer on the larger scale.  
Once we have this, we can iterate the argument
 to get estimates on larger and larger scales, eventually giving the strong rigidity.
Estimates proving  polynomial losses will be played off against estimates with exponential gains.

Four key ingredients in the proof of strong rigidity are:
\begin{enumerate}
\item   Gauge fixing.
\item New polynomial growth estimates for PDEs.
\item Propagation of almost splitting.
\item Quadratic rigidity in the right gauge. 
\end{enumerate}
We will use new polynomial growth estimates as ingredients in both (1) and (3).

\subsection{Gauge fixing}
Fix $(M,g,f)$. 
   We are given a diffeomorphism from a large compact set in $M$ to a second weighted space. 
\begin{itemize}
\item  The pull-back metric and weight are $g+h$ and $\e^{-f-k}$.
\item $h$ and $k$ are small on the compact set.
\end{itemize}
 Composing with a diffeomorphism on $M$ gives a different $h$ and $k$.   We want to mod out by this group action. 
  This is gauge fixing.

 One of the most interesting results of transformation groups is the existence of slices.
 A slice for the action of a group on a manifold  is  a submanifold  which is transverse to the orbits near a given point.   Ebin and Palais proved the existence of a slice for the diffeomorphism group of a {\emph{compact}} manifold acting on the space of all Riemannian metrics.  The slice can be thought of as the gauge fixing on the compact manifold.  
 However,  manifolds are not compact here. 

   In our setting, gauge fixing is choosing a diffeomorphism $\Phi$ on $M$ so the new $h$ is orthogonal to the group action.  
  This is a nonlinear PDE for $\Phi$.
   Orthogonality corresponds to making
 $\dv_f\,h=0$, where
  $\dv_f\,(h)=\e^{f}\,\dv\,(\e^{-f}\,h)=\dv\,(h)-h(\nabla f,\cdot)$.  
  The solution to the nonlinear PDE asks to find a diffeomorphism $\Phi$ so that   $\tilde{h} = \Phi^*\,(g+h)-g$ satisfies
\begin{align}
\dv_f\, \tilde{h} =0 \, .	\label{e:nonlinear}
\end{align}
  Terms involving $\dv_f $ come up again and again, so many quantities simplify in this gauge and having them drop out as they do when $\dv_f\, \tilde{h}=0$ makes things possible to analyze.

 We construct the diffeomorphism $\Phi$ that solves this nonlinear PDE using an iteration scheme  for the linearized operator $\cP$ on vector fields $Y$.    
 Using optimal polynomial  bounds on $\cP$, we show sharp polynomial bounds for the displacement function of $\Phi$
$$
x\to \text{dist}_g(x,\Phi (x))\,  .
$$
For applications, it is crucial that we only assume closeness on a compact set and, in particular, a priori 
the two shrinkers do not need to be diffeomorphic.  This means that we cannot fix the gauge at the outset.  Instead we need to apply our gauge fixing procedure iteratively to fix the gauge on larger and larger scales as we move outward and show closeness on larger and larger scales.  To pull this off requires very strong estimates for the displacement which is what we show.  Our optimal estimates show that the displacement of the gauge fixing diffeomorphism grows at a sharp polynomial rate.  The bound is relative so closeness on the initial scale implies closeness at a larger scale.  

On a shrinker $(M,g,f)$ there is a natural gaussian  $L^2$ norm given by  $\| u \|_{L^2}^2 = \int_M u^2 \, \e^{-f}$. 
Diffeomorphisms near the identity are infinitesimally generated by integrating a vector field $X$.  The infinitesimal change of the metric is given by the Lie derivative of the metric with respect to $X$.  This is equal to $- \frac{1}{2} \, \dv_f^* X$, where $\dv_f^*$ is the operator adjoint of $\dv_f$ with respect the to gaussian inner product. Thus,
if we define the  operator $\cP$ by
$$
\cP\,X=\dv_f\,\circ\dv_f^*\,X\,  ,
$$
then the linearization of \eqr{e:nonlinear} is to find a vector field $X$ with $$\cP\,X=\frac{1}{2}\,\dv_f\,h\,  .$$ 
A detailed analysis of $\cP$ and its properties play an important role in the gauge fixing.

The canonical second order elliptic ``drift Laplacian'' that is self-adjoint and positive definite with respect to the gaussian inner product is a generalized Ornstein-Uhlenbeck operator.   On tensors $u$ this operator is given by
$\cL \, u = \Delta \, u - \nabla_{\nabla f} u$.
  Given a vector field $X$ on a 
shrinker, the two operators $\cL$ and $\cP$ are related by the identity
$$
-2\, \cP\,X=\nabla \dv_f\,X+\cL\,X+\frac{1}{2}\,X\,  .
$$
The unweighted version of $\cP$  was used implicitly by Bochner to show that closed manifolds with negative Ricci curvature have no Killing fields. The unweighted operator was later used by Bochner and Yano, to show that the isometry group of such manifolds is finite.  
 The unweighted operator also arises in general relativity and fluid dynamics.    The weighted operator $\cP$ appears to have been largely overlooked.
The relationship between $\cP$ and the unweighted version, used by Bochner, mirrors the relationship between the Ornstein-Uhlenbeck operator and the Laplacian.  

\subsection{New polynomial growth estimates for PDEs}
Surprisingly, in very general settings, we show the same polynomial growth bounds for $\cL$ that Laplace and Hermite observed on   $\RR^n$  for the standard gaussian.    We use the relationship between $\cL$ and $\cP$ to translate these optimal bounds  for $\cL$ into optimal bounds for $\cP$.   These  optimal bounds   hold for all shrinkers for both Ricci and mean curvature flows, giving
a powerful  new tool for a wide array of problems.  These estimates are used in both the gauge fixing and the propagation of almost splitting.

\subsection{Propagation of almost splitting}
  One of the important new 
ingredients  is that a
   Ricci shrinker close to a product $N \times \RR^{n-\ell}$ on a large scale remains close on a fixed larger scale.  
The idea is that the initial closeness will imply that $\cL$ has  eigenvalues that are exponentially close to $\frac{1}{2}$. 
The drift Bochner formula on a shrinker implies that every eigenvalue is at least $\frac{1}{2}$ with equality only when it splits.  We show that being close to $\frac{1}{2}$ gives that the hessian is almost zero in $L^2$, which is very strong on the region where the weight $f$ is small but says almost nothing further out.  The crucial point is  that the hessian can grow only polynomially, so the very small initial bound
gives bounds much further out.  Thus, the gradients of these eigenfunctions give the desired almost parallel vector fields and almost splitting.
  

If a shrinker is exponentially close to a cylinder on scale $R$, then almost splitting gives
on scale $(1+ \epsilon) \, R$:
 \begin{enumerate}
\item[(A)] $n-\ell$ almost translations and a metric almost splitting.

\item[(B)] The slices $\{ f = c\}$ are almost spherical.

\item[(C)] $f$ also almost splits $f=f_0 + \frac{|x|^2}{4}$.

\end{enumerate}
Combining (A) and (B)  gives  a diffeomorphism on  scale $(1+ \epsilon) \, R$ from the cylinder to the second shrinker and this diffeomorphism is close to an isometry.  The condition (C) guarantees that the diffeomorphism also almost preserves $f$.
As a result, we see that the shrinker looks cylindrical even on the larger scale.  However, there is a loss in the estimates - it may look less 
cylindrical on the larger scale - that makes this impossible to iterate on its own.
   
\subsection{Quadratic rigidity}

The propagation of almost splitting and gauge fixing give that the shrinker is  close to a cylinder on a large set via a diffeomorphism  that fixes the gauge.  
The last of the four key ingredients is an estimate for the difference in metrics that is small enough to be iterated.
For this, it is essential   that the gauge be right, or else it just isn't true.  The closeness cannot be seen via
 linear analysis.  However, we show that there is a second order rigidity that gives the estimate; we call this quadratic rigidity.

To explain the estimate,  let $(M,g,f)$ be the cylinder and  $(M,g+h, f+k)$ the shrinker that is close on a large compact set.  We need
 bounds on $h$ and $k$ that can be iterated.
The linearization of the shrinker equation is
\begin{align}	\label{e:linea}
	\frac{1}{2} \, L \, h   + \Hess_{ \frac{1}{2} \, \Tr (h) - k} + \dv_f^* \, \dv_f \, h \, .
\end{align}
This linearization was derived by Cao-Hamilton-Ilmanen in their calculation of 
 the second variation operator for Perelman's  entropy.
The operator $L$ acts on $2$-tensors by $$L \, h = \cL \, h + 2 \, \R (h)\, ,$$
and $\R(h)$ is the natural action of the Riemann tensor.

Since $(M,g+h, f+k)$ is also a shrinker,  \eqr{e:linea} must be at least quadratic in $(h,k)$.  
The last two terms in \eqr{e:linea} are gauge terms -  i.e., in the image of $\dv_f^*$ and there is no reason for these - or $h$ -  to be small if not in right gauge.
 In the {\it right gauge}, the difference $h$ in the metrics satisfies the Jacobi equation $L\, h = 0$ up to higher order terms.  
 This does not  force $h$ to be small since
 cylinders have non-trivial Jacobi fields that could potentially integrate to give nearby shrinkers.  
The second variation of the {\bf{shrinker equation}} in the direction of a Jacobi field is given by the tensor
    \begin{align}	\label{e:quad}
    	-2\, |\nabla u|^2 \, \Ric - 2\, S \, u \, \Hess_u - S \, \nabla u \otimes \nabla u
	\, . 
\end{align}
Here  $S$ is scalar curvature and $u$ is a quadratic Hermite polynomial that
 measures the projection of $h$ onto Jacobi fields.  
 The first order Taylor expansion will give that $h$ is a Jacobi field to first order and, thus, $|h|$ is  $|u|$ up to higher order.  On the other hand, the second order Taylor expansion
 will imply that 
\eqr{e:quad} 
vanishes to at least third order in $h$.  Combining these, we see that the quadratic expression \eqr{e:quad} is in fact at least
 cubic in $u$.  When $u$ is small, this implies that $u$ and $h$ vanish; we will have extra error terms  so 
 will get that $h$ is exponentially small, giving the improvement that we needed to iterate.

\end{document}